\documentclass[titlepage,a4paper]{article}
\usepackage[T1]{fontenc}
\usepackage[latin1]{inputenc}
\usepackage{amsmath}
\usepackage{amssymb}

\usepackage{amsthm}

\newtheorem{ex}{Example}
\newtheorem{Def}{Definition}

\newtheorem{rem}{Remark}
\newtheorem{Th}{Theorem}

\newtheorem{lem}{Lemma}

\newtheorem{cor}{Corollary}

\newcommand{\RR}{\mathbb{R}}

\newcommand{\NN}{\mathbb{N}}

\newcommand{\cK}{{\cal K}}
\newcommand{\cU}{{\cal U}}

\newcommand{\cR}{{\cal R}}

\newcommand{\cI}{{\cal I}}

\newcommand{\cD}{{\cal D}}

\newcommand{{{\cadlag}}}{c\`adl\`ag}

\newcommand{\id}{\mathrm{id}}

\begin{document}

\sffamily
\sloppy
\pagestyle{plain}

\renewcommand{\baselinestretch}{1}

\title{Approximate Bermudan option pricing based on the r\'eduite or cubature: soundness and characterisation of perpetual prices as fixed points}
\author{Frederik S Herzberg\thanks{Abteilung f\"ur Stochastik, Institut f\"ur Angewandte Mathematik, Universit\"at Bonn, D-53115 Bonn, Germany ({\tt herzberg@wiener.iam.uni-bonn.de})} \thanks{Mathematical Institute, University of Oxford, Oxford OX1 3LB, England} }
\date{}

\maketitle

\renewcommand{\baselinestretch}{1}

\begin{abstract} 

In this paper, it is shown that Bermudan option pricing based on either the r\'eduite (in a one-dimensional setting: piecewise harmonic interpolation) or cubature -- is sensible from an economic vantage point: Any sequence of thus-computed prices for Bermudan options with increasing sets of exercise times is increasing. Furthermore, under certain regularity assumptions on the payoff function and provided the exercise times are equidistant of exercise mesh size $h$, it has a supremum which coincides with the least fixed point of the approximate pricing algorithm -- this algorithm being perceived as a map that assigns to any real-valued function $f$ (on the basket of underlyings) the approximate value of the European option of maturity $h$ and payoff function $f$.

\end{abstract}

\noindent

\section{Introduction}

The Bermudan pricing approach outlined in the first part of this paper is to iteratively construct a piecewise {\em harmonic} approximation to the function that assigns the value of a Bermudan option (with payoff function $g$, exercise mesh size $t>0$ and maturity $nt$) to the vector of logarithmic start prices of the underlying assets. 

In the first step of this iteration, one will compute a piecewise harmonic approximation to the function that assigns the corresponding European option price of maturity $t$ to the logarithmic asset prices at the penultimate time $T-t$ where exercise is possible. The iteration step will consist in computing the expectation of this function after time $t$ (under the assumption that the vector price process is a Markov process), discount, take the maximum with the payoff function $g$, and then perform a r\'eduite-based interpolation (in the one-dimensional setting: an interpolation by interval-wise harmonic functions). (This procedure computes a piecewise harmonic approximation to the so-called {\em Snell envelope} \cite{Snell} from the theory of optimal stopping \cite{Nev,ElK}.)

The notion of harmonicity will be derived from the log-price process $X$, assumed to be a time-continuous Markov process with a path-continuous modification, i e a diffusion -- which then makes the infinitesimal generator of $X$ a second-order elliptic differential operator and gives rise to a notion of harmonicity.

The principal advantage of such an algorithm based on piecewise harmonic interpolation is that it is intrinsically finite-dimensional, since the space of harmonic functions corresponding to $X$ will be two-dimensional. However, one needs to make sure that this algorithm will make sense economically. In particular one will have to require that the iteration step be monotone, since adding one more possible exercise time of course increases the value of the option price. This is not difficult to see and will be one of the first results of this paper. Also, once monotonicity of the iteration step has been established, the question of whether there exists a supremum to this sequence of approximate non-perpetual option prices arises, since this could then be treated as an approximation to the price of the corresponding perpetual option of exercise mesh size $t$. Under certain assumptions on the payoff function, this will be proven as well; furthermore we will be able to characterise it as the minimal fixed point of the iteration step previously referred to.

First, we will discuss these questions in the one-dimensional setting -- very little knowledge of potential theory has to be assumed for the proofs in that section. Second, we shall generalise that approach to higher dimensions; this will entail a few technical subtleties. At the end of both Sections we will apply these results to financially relevant settings.

Finally, in a last Section we shall discuss related questions for Bermudan pricing based on cubature and establish a linear convergence rate that corresponds to the discount factor.

\section{Some terminology for Bermudan option pricing algorithms}

Let $C^0\left(\RR^d,[0,+\infty)\right)$ and $L^0\left(\RR^d,[0,+\infty)\right)$, as usual, denote the spaces of nonnegative continuous functions defined on $\RR^d$, and of equivalence classes of nonnegative measurable functions defined on $\RR^d$, respectively.

The observed price for any financial derivative increases when adding one more possible exercise time; also, the perpetual limit often has to be finite (in particular for puts, since then any option price is bounded by the strike price). 

Approximative Bermudan pricing algorithms also should reflect these features of actual observed prices, therefore we introduce the following manner of speaking to distinguish economically sensible from less sensible algorithms.

\begin{Def}\label{algorithmclasses} A map $D:L^0\left(\RR^d,[0,+\infty)\right)\rightarrow L^0\left(\RR^d,[0,+\infty)\right)$ is said to be a {\em sound iterative Bermudan option pricing algorithm} (for short, a {\em sound algorithm}) for a payoff function $g\in C^0\left(\RR^d,[0,+\infty)\right)$ if and only if $Df\geq g$ for all $f\in L^0\left(\RR^d,[0,+\infty)\right)$ and the map $D$ is pointwise monotone, that is \begin{eqnarray*}&& \forall f_0,f_1\in L^0\left(\RR^d,[0,+\infty)\right)\\ && \left(\forall x\in\RR^d \quad f_0(x)\leq f_1(x)\Rightarrow \forall x\in\RR^d \quad Df_0(x)\leq Df_1(x) \right).\end{eqnarray*} A sound iterative Bermudan option pricing algorithm $D$ is said to {\em have a perpetual limit} if and only if $$\sup_{n\in\NN_0} D^{\circ n}g\in L^0\left(\RR^d,[0,+\infty)\right)$$ (rather than this supremum being allowed to equal $+\infty$ on a subset of positive Lebesgue measure of its range). In that very case, the function in the last line is simply referred to as the {\em perpetual limit} of the algorithm. Finally, $D$ is said to {\em converge linearly in $L^\infty$ to the perpetual limit} if and only if there exists a $c\in(0,1)$ such that $$\forall n\in\NN\quad \left\|\left(D^{n+1}-D^n\right)g\right\|_{L^\infty(\RR^d,\RR)}\leq c\cdot \left\|\left(D^{n}-D^{n-1}\right)g\right\|_{L^\infty(\RR^d,\RR)},$$ $D^n$ being shorthand for $D^{\circ n}$ for all $n\in\NN_0$.
\end{Def}

\begin{rem} The elements of $L^0\left(\RR^d,[0,+\infty)\right)$ should be conceived of as assigning the value -- that is, the expected payoff -- of an option to the vector of logarithmic start prices of the components of the basket (at least on the complement of a Lebesgue null set).
\end{rem}

\begin{rem} \label{soundmonotone} The monotonicity condition imposed on sound iterative Bermudan pricing algorithms entail that the sequence of functions $\left(D^ng\right)_{n\in\NN_0}=\left(D^{\circ n}g\right)_{n\in\NN_0}$ is always pointwise increasing. Thus, this sequence has a limit: $$g\leq \sup_{n\in\NN} D^ng= \lim_{n\rightarrow\infty} D^ng.$$
\end{rem}

The infimum of all $D$-fixed points is always an upper bound for the perpetual limit: 

\begin{lem} \label{inffixedpoint}Let $D$ be a sound iterative Bermudan pricing algorithm for $g$ with a perpetual limit $u$. Then the function $u$ is smaller than any fixed point of $D$; moreover, $Du\geq u$.
\end{lem}
\begin{proof} Any fixed point $h$ of $D$ is in the image of $D$ and therefore, due to our assumptions on sound algorithms, pointwise greater than or equal to $g$. Now, as $D$ (and thus $D^n$) is pointwise monotone, $$\forall n\in\NN_0 \quad h=D^nh\geq D^ng,$$ therefore $$h=\sup_m D^mh\geq \sup_mD^mg,$$ where the right hand side is just the perpetual limit. Hence, any fixed point of $D$ is greater than or equal to the perpetual limit. Furthermore, observe that due to the pointwise monotonicity of $D$, $$D\left(\sup_nD^ng\right)\geq D^{m+1}g,$$ therefore for all $m\in\NN_0$, $$D\left(\sup_nD^ng\right)\geq \sup_m D^{m+1}g=\sup_mD^mg.$$ 
\end{proof}

Later on, it will turn out that if $D$ is based on piecewise harmonic interpolation or the r\'eduite, the perpetual limit is, in fact the minimal fixed point (cf Theorem \ref{fixexist} and Lemma \ref{Lemma7.6} for piecewise harmonic interpolation, Theorem \ref{Theorem7.2} for r\'eduite-based approximation, and Theorem \ref{monoconv} for a result on a map $D$ which is based on cubature).

\section{Piecewise harmonic Bermudan option pricing for options on one asset}

\subsection{Definitions and facts from classical potential theory}

Let $(P_t)_{t\geq 0}$ be a Markovian semigroup associated to a one-dimensional diffusion and let $L$ be the infinitesimal generator of $L$. Then $L$ will be a second-order differential operator (cf eg Revuz and Yor \cite{RY}).

If one assumed, for the sake of an example, that $L$ is a {second-order differential operator} even with constant coefficients, that is $$L:f\mapsto \beta f' + \frac{1}{2} f'',$$ then $P$ will merely be the semigroup of finite-dimensional distributions of a multiple of Brownian motion with linear drift.

Let now $U\subseteq \RR$ be open and nonempty. A twice differentiable function $f:U\rightarrow\RR$ is called {\em harmonic} on $U$ if and only if $Lf=0$ on $U$. A continuous function $g:\bar U\rightarrow \RR$ is said to be {\em subharmonic on }$U$ if and only if $Lg$ is defined as an element of $\RR\cup\{\pm\infty\}$ and $Lg(x)=\lim_{t\downarrow 0}\frac{P_tg(x)-g(x)}{t}\geq 0$ for all $x\in U$. A function $f:\bar U\rightarrow \RR$ will be called {\em superharmonic on }$U$ if and only if $-f$ is subharmonic on $U$.

A function that is harmonic (subharmonic, superharmonic) on $\RR$ will simply be called {\em harmonic} ({\em subharmonic, superharmonic}).

Also, if $h$ is both superharmonic and subharmonic, then one will have $h\geq P_th\geq h$ for all $t\geq 0$, making $h$ harmonic.
 
The space of all harmonic functions on $U$ is a vector space. As a classical result (cf eg Protter and Weinberger \cite{PW}), one has that all subharmonic functions obey the following

{\em Maximum Principle: If $g$ is subharmonic on $U\subset \RR$ and $\bar U$ is compact, then the maximum of $g$ will be attained on $\partial U$.}

Also, the set of subharmonic functions has the following closure properties:

{\em Subharmonic functions form a cone, closed under $\vee$ and $P$: If $f,g:\RR\rightarrow \RR$ are subharmonic functions, then so are $f\vee g$, $P_tf$ for all $t\geq 0$, and $\alpha f+\beta g$ for all $\alpha,\beta\geq 0$. }

See textbooks like Revuz and Yor \cite{RY}, Meyer \cite{M}, Port and Stone \cite{PS} and the classical references therein.

\subsection{Harmonic interpolation}

\begin{lem} \label{uniqueinterpol} Given two support abscissas $a_0 \neq a_1$ and ordinates $c_0$, $c_1$, there is a unique harmonic interpolation, that is a harmonic function $h:\RR\rightarrow\RR$ such that $h(a_i)=c_i$ for $i\in\{0,1\}$. Thus, there cannot be more than two linearly independent harmonic functions, and two harmonic functions coincide on all of $\RR$, once they agree in two points. 

Assume now, for the sake of simplicity, that $L$ is a {second-order differential operator}  with constant coefficients, that is $$L:f\mapsto \beta f' + \frac{1}{2} f''.$$ Then, if $\beta\neq 0$, all harmonic functions have the shape $\gamma_0 + \gamma_1\cdot \exp\left(-2\beta\cdot \right)$ for some real constants $\gamma_0,\gamma_1$, and if $\beta =0$, the vector space of all harmonic functions is just the space of all affine functions.
\end{lem}
\begin{proof} The existence and uniqueness assertion for harmonic interpolation is merely the statement of the unique solvability of the Dirichlet problem for $L$. 

If there were more than two linearly independent harmonic functions on $(a_0,a_1)$, the operation of harmonic interpolation would not be unique. 

Finally, constant functions are of course harmonic. Also, in case $\beta =0$ any linear function must be harmonic as well, whence $\id_\RR$ and $1$ are already to linearly independent functions for the case $\beta=0$. For the case of non-zero $\beta$, we remark that $$L\exp\left(-2\beta \cdot\right) = -2\beta\exp\left(-2\beta \cdot\right)\cdot\beta + 4\beta^2\exp\left(-2\beta \cdot\right)\cdot \frac{1}{2} =0 $$ and conclude that $1$ and $\exp\left(-2\beta \cdot\right)$ are two (obviously linearly independent) harmonic functions.

\end{proof}

\begin{lem} \label{threezeroes} Let $I=(a,b)$ be an open nonempty interval, and let $f:\bar I\rightarrow\RR$ be subharmonic on $U$ such that $f(a)=f(b)=0$. If there exists an $x_0\in I$ with $f(x_0)=0$, then $f=0$ on $I$.
\end{lem}
\begin{proof} Due to the Maximum Principle, $f\leq 0$. Now suppose $f(x_0)=0$ for some $x_0\in(a,b)$. If neither $f=0$ on $(a,x_0)$ nor $f=0$ on $(x_0,b)$, then we would have a $y_0\in (a,x_0)$ and a $y_1\in (x_0,b)$ such that $f(y_0)<0$ and $f(y_1)<0$. Then $\max_{\partial(y_0,y_1)}f<\max_{[y_0,y_1]} f =f(x_0)=0$, a contradiction to the Maximum Principle. Hence, either $f=0$ on $(a,x_0)$, or $f=0$ on $(x_0,b)$, and without loss of generality, we shall from now on assume $f=0$ on $(a,x_0)$.

If there was now an $x_1\in (x_0,b)$ with $f(x_1)<0$, pick any $a_1\in (a,x_0)$ and let (by virtue of Lemma \ref{uniqueinterpol}) $h_1:\RR\rightarrow\RR$ be the unique harmonic function such that $h_1(a_1)=f(a_1)=0$, $h_1(x_1)=f(x_1)<0$. Then $f-h_1$ is subharmonic on $(a_1,x_1)$ and vanishes on the boundary of $(a_1,x_1)$. Hence, by the Maximum Principle applied to $f-h_1$, one has $f-h_1\leq 0$ on $(a_1,x_1)$. However, $f=0$ on $(a_1,x_0)\subset (a,x_0)$ as established previously, and, by the Maximum Principle applied to $h_1$, also $h_1\leq 0$ on $(a_1,x_1)\supset (a_1,x_0)$. Thus $$0\leq -h_1 = f-h_1\leq 0 \text{ on }(a_1,x_0),$$ which means $h_1=0$ on $(a_1,x_0)$. But since $(a_1,x_0)$ is a nonempty open interval (and therefore does not have less than two elements) and $h_1$ is harmonic, this can only be true if $h_1$ is identically zero on all of $\RR$, in contradiction to $h_1(x_1)<0$.

\end{proof}

If there is only a first order non-zero term, the space of harmonic functions will just coincide with the space of constant functions.

\begin{lem}\label{subinterpol} 

Let $a_0<\dots<a_m\in\RR$. Consider a function $f:\RR\rightarrow\RR$ that is subharmonic on $(a_0,a_m)$ and let, by virtue of Lemma \ref{uniqueinterpol}, $\cI(f)$ be the unique interpolation to $f$ such that $\cI(f)=f$ on $\{a_0,\dots,a_m\}$ and for all $i<m$, $\cI(f)$ coincides with some harmonic function on $(a_{i},a_{i+1})$, denoted by $f_i$. 

Then the following assertions hold:

\begin{enumerate}
\item $\cI(f)\geq f$ on $[a_0,a_m]$.

\item $\cI(f)=\max\{f_0,\dots,f_{m-1}\}$ on $(a_0,a_m)$, thus making $\cI(f)$ subharmonic on $(a_0,a_m)$.

\item $\cI(f)\leq f$ on $(-\infty,a_0)\cup (a_m,+\infty)$. 

\end{enumerate}

\end{lem}

\begin{cor} \label{interpolsubharmonic}Piecewise harmonic interpolation with respect to a set of support abscissas $\{a_0,\dots,a_m\}$ preserves subharmonicity on $(a_0,a_m)$ in the following sense: If, in the notation of Lemma \ref{subinterpol}, $f$ is subharmonic on $(a_0,a_m)$, then so is $\cI(f)$.
\end{cor}

\begin{proof}[Proof of Lemma \ref{subinterpol}]
\begin{enumerate} 
\item We apply for every $i<m$ the Maximum Principle to the subharmonic function $f-f_i$ on $(a_i,a_{i+1})$ (which vanishes on the boundary of $(a_i,a_{i+1})$) to get $0\geq f-f_i =f-\cI(f)$ on $[a_i,a_{i+1}]$ for all $i<m$. This proves $\cI(f)\geq f$ on $[a_i,a_{i+1}]$. 

\item Again from the Maximum Principle, we get that two different harmonic functions can only have one point in common. Hence, not only $\left\{f_i=f_{i+1}\right\}\supseteq \left\{a_{i+1}\right\}$ for all $i<m-1$ (a trivial consequence of the interpolation) but if $f_i\neq f_{i+1}$, then even $$\left\{f_i=f_{i+1}\right\}=\left\{a_{i+1}\right\}.$$ Consider any $i<m$. There are two possibilities: Either $f_i< f_{i+1}$ on $(-\infty,a_{i+1})$ and $f_i> f_{i+1}$ on $(a_{i+1},+\infty)$ or the other way round $f_i> f_{i+1}$ on $(-\infty,a_{i+1})$ and $f_i< f_{i+1}$ on $(a_{i+1},+\infty)$. 

Suppose, for a contradiction, the former situation holds for some $i<m$: $f_i< f_{i+1}$ on $(-\infty,a_{i+1})$ and $f_i> f_{i+1}$ on $(a_{i+1},+\infty)$. Then $\cI(f)$ would equal $f_i\wedge f_{i+1}$ on $[a_i,a_{i+2}]$, which is superharmonic. Then, $f-\left(f_i\wedge f_{i+1}\right)$ would be subharmonic on $[a_i,a_{i+2}]$ and it would have three zeroes, in $a_i$, $a_{i+1}$ and $a_{i+2}$. By Lemma \ref{threezeroes}, this can only be true if $f-\left(f_i\wedge f_{i+1}\right)=0$ on all of $[a_i,a_{i+2}]$. Thus, $f=f_i\wedge f_{i+1}$ on $[a_i,a_{i+2}]$. Since $f$ is subharmonic on $[a_i,a_{i+2}]$, so must then be $f_i\wedge f_{i+1}$ (which is already superharmonic on $[a_i,a_{i+2}]$) then, and therefore, $f_i\wedge f_{i+1}$ is harmonic on $[a_i,a_{i+2}]$. Now, $f_i$ and $f_{i+1}$ agree only on one point, viz. $a_{i+1}$, inside $(a_i,a_{i+2})$. Therefore, $f_i\wedge f_{i+1}$ has with both $f_i$ and $f_{i+1}$ at least two points in common. Since any harmonic function is uniquely determined by a mere two points, we thus get $$f_i =f_i\wedge f_{i+1}=f_{i+1}$$ in contradiction to our assumption.
 
Therefore, $f_i> f_{i+1}$ on $(-\infty,a_{i+1})$ and $f_i< f_{i+1}$ on $(a_{i+1},+\infty)$ for all $i<m$. 

Thus, due to $a_0<\dots<a_m$, we gain \begin{align*}f_0&>f_1 &\text{ on } & (-\infty,a_{1})\supset(-\infty,a_{0}),\\ f_1&>f_2 &\text{ on }& (-\infty,a_{2})\supset\dots\supset(-\infty,a_{0}), \\ &\vdots &\vdots & \vdots \\ f_{m-2}&>f_{m-1}& \text{ on } &(-\infty,a_{m-1})\dots\supset(-\infty,a_{0}) \end{align*} as well as \begin{align*}f_0&<f_1 &\text{ on }& (a_{1}, +\infty)\supset\dots\supset(a_{m},\infty),\\ f_1&<f_2 &\text{ on } &(a_{2}, +\infty)\supset\dots\supset(a_{m},\infty), \\ &\vdots &\vdots &\vdots \\ f_{m-2}&<f_{m-1} & \text{ on } &(a_{m-1}, +\infty). \end{align*} Hence we arrive at \begin{align*}&f_0>f_1>\dots>f_{m-1} \text{ on } (a_{0},a_1)\end{align*} as well as 
\begin{align*}\forall i\in\{1\dots,m-1\}\quad \left(f_0<\dots <f_i, \quad f_i>\dots>f_{m-1}\right) \text{ on } (a_{i},a_{i+1})\end{align*} Therefore, $$ \forall i<m\quad \forall j\neq i\qquad  f_i>f_j\text{ on } (a_{i},a_{i+1}),$$ so $$ \forall i<m\quad f_i=\max\left\{f_0,\dots,f_{m-1}\right\} \text{ on } (a_{i},a_{i+1})$$ which proves $\cI(f)=\max\left\{f_0,\dots,f_{m-1}\right\}$.

\item This part of the proof is similar to the proof of Lemma \ref{threezeroes}. The function $f-\cI(f)$ is subharmonic on $(-\infty,a_1)$ and it has two zeroes in $a_0$ and $a_1$. First consider the case where $f$, whilst being subharmonic, is not harmonic on $(a_0,a_1)$, which ensures that the subharmonic function $f-\cI(f)$ -- which vanishes on the boundary of $(a_0,a_1)$ -- is not harmonic. Then there is an $x_1\in (a_0,a_1)$ such that $\left(f-\cI(f)\right)(x_1)<0$. If $f-\cI(f)$ was now also negative for some $x_0\in (-\infty,a_0)$, then the Maximum Principle would yield $0>\max_{\{x_0,x_1\}}\left(f-\cI(f)\right)=\max_{[x_0,x_1]}\left(f-\cI(f)\right)$, but on the other hand $\left(f-\cI(f)\right)(a_0)=0$ and $a_0\in (x_0,x_1)$, a contradiction. Hence in the case where $f$ is not harmonic on $(a_0,a_1)$, we already have $f-\cI(f)\geq 0$ on $(-\infty,a_0)$. 

Next suppose $f$ is indeed harmonic on $(a_0,a_1)$. Then, so will be $f-\cI(f)$ and because $f-\cI(f)=0$ on $\{a_0,a_1\}$, we will get $f-\cI(f)=0$ on all of $[a_0,a_1]$. Now pick some $x_1\in (a_0,a_1)$ and assume, for a contradiction that $\left(f-\cI(f)\right)(x_0)<0$ for some $x_0\in (-\infty,a_0)$. According to Lemma \ref{uniqueinterpol}, there exists a unique harmonic function $h$ such that $h(x_0)=\left(f-\cI(f)\right)(x_0)<0 $ and $h(x_1)=\left(f-\cI(f)\right)(x_1)=0$. Then $f-\cI(f)-h$ is subharmonic on $(x_0,x_1)$ and vanishes on the boundary of $(x_0,x_1)$, hence $f-\cI(f)-h\leq 0$ on $[x_0,x_1]$. But also, $h$ is harmonic and nonpositive on the boundary of $(x_0,x_1)$, thus $-h$ will be nonnegative on all of $(x_0,x_1)$. Therefore $$0\leq -h = f-\cI(f)-h\leq 0 \text{ on }[a_0,x_1],$$ therefore $h=0$ on $(a_0,x_1)\neq \emptyset$, yielding $h=0$ everywhere. This contradicts $h(x_0)<0$ and hence completes the proof for $f-\cI(f)\geq 0$ on $(-\infty,a_0)$.

Analogously, one can prove the domination of $\cI(f)$ by $f$ on $(a_m,+\infty)$.

\end{enumerate}
\end{proof}

\begin{lem}\label{monointerpol} Piecewise harmonic interpolation to a set of support absicssas $\{a_0,\dots,a_m\}$ is monotone on $[a_0,a_m]$ in the sense that if $f\leq g$ on $[a_0,a_m]$, then the piecewise harmonic interpolation $\cI(f)$ of $f$ will be dominated by the piecewise harmonic interpolation $\cI(g)$ of $g$ on $[a_0,a_m]$.
\end{lem}
\begin{proof} Consider any $i<m$. By defintion, $\cI(f)=f\leq g=\cI(g)$ on $\{a_i,a_{i+1}\}$, hence $\cI(f)-\cI(g)\leq 0$ on $\{a_i,a_{i+1}\}$. But $\cI(f)-\cI(g)$ is harmonic on $(a_i,a_{i+1})$, therefore the Maximum Principle yields $\cI(f)-\cI(g)\leq 0$ on $[a_i,a_{i+1}]$. Since this holds for all $i<m$, we get $\cI(f)-\cI(g)\leq 0$ on all of $[a_0,a_m]$
\end{proof}

\begin{lem}\label{monointerpolsub} Let, as before, $\cI:\RR^{[a_0,a_m]}\rightarrow \RR^\RR$ denote the operator of piecewise harmonic interpolation with respect to the set of support abscissas $\{a_0,\dots,a_m\}$. Consider a subharmonic function $f:\RR\rightarrow\RR$ and a harmonic function $h:\RR\rightarrow\RR$ such that $f\leq h$ on $\RR$. Then $\cI(f)\leq h$ on $\RR$.
\end{lem}
\begin{proof} From the previous Lemma \ref{monointerpol}, we already know that $\cI(f)(x)\leq\cI(h)(x)$ holds for all $x\in [a_0,a_m]$. However, $\cI(h)=h$, hence $\cI(f)\leq h$ on $[a_0,a_m]$ and from Lemma \ref{subinterpol}, we conclude that $h\geq f\geq\cI(f)$ on the intervals $(-\infty,a_0)$ and $(a_m,+\infty)$.
\end{proof}

\subsection{A fixed point theorem}

\begin{Th}\label{fixexist} Let $\cI:\RR^{[a_0,a_m]}\rightarrow \RR^\RR$ again denote the operator of piecewise harmonic interpolation with respect to the set of support abscissas $\{a_0,\dots,a_m\}$, and fix $t>0$, $r\geq 0$. Let $g$ and $c$ be subharmonic functions, and let $h$ be harmonic and nonnegative. Let $c$ be, moreover, harmonic on each of the intervals $[a_i,a_{i+1}]$ for $i<m$, and assume $P_tc(x)$ is finite for all $x\in\RR$. Suppose $c\leq g$ on $[a_0,a_m]$ and $c,g\leq h$ on $\RR$. Now define $$\cK:f\mapsto\cI\left(e^{-rt}P_t\left(\cI(f)\vee c\right)\vee g\right)\restriction [a_0,a_m]$$ as well as $$Q:=\left\{f\restriction [a_0,a_m]\ : \ \begin{array}{c} f:\RR\rightarrow\RR \text{ subharmonic}, \quad f\geq c \text{ on }[a_0,a_m], \\ \forall i\in\{1,\dots,m-2\} \quad f\text{ harmonic on }(a_i,a_{i+1}), \\ f\text{ harmonic on }(-\infty,a_1),(a_{m-1},+\infty), \quad f\leq h \end{array}\right\}.$$ Then $\cK$ maps the convex and bounded subset $Q$ of $C^0[a_0,a_m]$ continuously to itself. Moreover, due to Lemma \ref{uniqueinterpol}, $Q$ is a subset of a finite-dimensional subspace of $C^0[a_0,a_m]$ (this subspace being the space of all functions from $[a_0,a_m]$ that are harmonic on each of the intervals $[a_i,a_{i+1}]$ for $i<m$. By Brouwer's Fixed Point Theorem, $\cK$ has got a fixed point in $Q$. Finally, $\cK$ is a composition of monotone functions on $[a_0,a_m]$ and therefore monotone as well.
\end{Th}
\begin{proof}[Proof] We can divide the proof for $\cK(Q)\subseteq Q$ into three parts:
\begin{enumerate}
\item The cone of subharmonic functions is closed under $\vee$, under $P_t$, under multiplication by constants and under piecewise harmonic interpolation $\cI$ (cf Lemma \ref{subinterpol}), therefore the image of $Q$ under $\cK$ can only consist of subharmonic functions. 
\item By the monotonicity of $P_t$ and $\cI$ (Lemma \ref{monointerpol}), combined with the equations $P_th=h$ and $\cI(h)=h$, as well as $c, g\leq h$, we have for all $f:\RR\rightarrow\RR$ such that $f\restriction [a_0,a_m]\in Q$, $$ e^{-rt}P_t\left(\underbrace{\cI(h)\vee c}_{=h\vee c = h}\right)\vee g =  e^{-rt}\underbrace{P_th}_{=h\geq 0}\vee g\leq h.$$ But the left hand side of this last estimate is subharmonic (because the cone of subharmonic functions is closed under $\vee$, $P_t$ and, as we have remarked in Corollary \ref{interpolsubharmonic}, also $\cI$). Hence, due to Lemma \ref{monointerpolsub}, we obtain $\cK f\leq h$.
\item The lower bound follows again from the monotonicity of $\cI$ (Lemma \ref{monointerpol}, that is), but this time only by exploiting $c\leq g$ on $[a_0,a_m]$ and the specific shape of $c$: For all $f:\RR\rightarrow\RR$ such that $f\restriction [a_0,a_m]\in Q$, one has $$\cK f = \cI(\dots \vee g)\geq  \cI(g) \geq \cI(c)=c.$$ 
\end{enumerate}
We get that $Q$ is bounded by $\sup_{[a_0,a_m]}h \vee |c|\geq 0$ as a subset of $C^0[a_0,a_m]$, and because $Q$ is finite-dimensional, we may apply Schauder's Theorem, as soon as we have established the continuity of $\cK$ with respect to the norm $\|\cdot \|_Q:f\mapsto \max{\left\{\left|f(a_0)\right|,\dots,\left|f(a_m)\right|\right\}}$. 

Note, for this sake, that from the Maximum Principle, we have $\left\|f\right\|_Q=\max_{[a_0,a_m]}|f|$ whenever $f$ is harmonic on each of the intervals $(a_i,a_{i+1}$ for $i<m$. Therefore, if $(f_n)_{n\in\NN}\in Q^\NN$, $f\in Q$ and $\left\|f_n-f\right\|_Q\longrightarrow 0$ for $n\rightarrow \infty$, then $f_n\longrightarrow f$ uniformly on $[a_0,a_m]$ as $n\rightarrow \infty$, hence -- since all $f_n$ and $f$ are harmonic on $(-\infty,a_1)$ and $(a_{m-1},+\infty)$ -- by Harnack's inequality (cf eg \cite[p 111f]{PW}), $f_n\longrightarrow f$ pointwise on all of $\RR$ as $n\rightarrow \infty$. But $$c(x) \leq f_n(x) \leq h(x)=P_th(x)$$ for all $n\in\NN$ and $x\in\RR$, therefore Lebesgue's Dominated Convergence Theorem (applied to the measure $A\mapsto P_t\chi_A(x)$ for each $x\in\RR$) may be employed to get $P_tf_n(x)\longrightarrow P_tf(x)$ as $n\rightarrow\infty$ for all $x\in\RR$. Via \begin{align*}&\left|e^{-rt}P_t\left(\cI(f)\vee c\right)\vee g- e^{-rt}P_t\left(\cI(f_n)\vee c\right)\vee g\right|\\ \leq & \left|e^{-rt}P_t\left(\cI(f)\vee c\right)- e^{-rt}P_t\left(\cI(f_n)\vee c\right)\right|\\ \leq & \left|e^{-rt}P_t\left(\cI(f)\vee c- \cI(f_n)\vee c\right)\right| \\ \leq & \left|e^{-rt}P_t\left(\cI(f)-\cI(f_n)\right)\right| \\ = & \left|e^{-rt}P_t\left(f-f_n\right)\right| ,\end{align*} we can now deduce that $\left|\cK f_n- \cK f\right\|_Q \longrightarrow 0$ as $n\rightarrow\infty$. 

\end{proof}

The existence of a minimal fixed point for $\cK$ can be proven constructively as well:

\begin{cor} Let us adopt the notation of the previous Theorem. Then the sequence $\left(\cK^n(g\vee 0)\right)_{n\in\NN_0}$ is monotone on $[a_0,a_m]$, bounded and dominated by $h$. Therefore we have the existence of a limit on $[a_0,a_m]$ given by $$\forall x\in[a_0,a_m]\quad q(x):=\lim_{n\rightarrow \infty}\cK^n(g\vee 0)(x)=\sup_{n\in\NN_0}\cK^n(g\vee 0)(x).$$ This limit is an element of $Q$ and therefore can be canonically extended to the whole of $\RR$. By the continuity of $\cK$, $q$ is a fixed point of $\cK$. On $[a_0,a_m]$, the convergence in the last equation will be uniform.
\end{cor}
\begin{proof} The only part of the Corollary that does not follow directly from the preceding Theorem \ref{fixexist} is the uniformity of the convergence and that $q$ will be harmonic on each of the intervals $[a_i,a_{i+1}]$ for $i<m$. However, monotone convergence on compact sets preserves harmonicity and is always uniform (cf e g Meyer \cite{M} -- or, more directly, Port and Stone \cite[Theorem 3.9]{PS} if $P$ is the Brownian semigroup).
\end{proof}

\begin{lem}\label{Lemma7.6}In the preceding Corollary's notation, $q$ is the minimal nonnegative fixed point of $\cK$.
\end{lem}
\begin{proof} The proof partly copies the one for Lemma \ref{soundmonotone}. Any nonnegative fixed point $p$ of $\cK$ must be greater or equal $g$ on $[a_0,a_m]$. Therefore the monotonicity of $\cK$ on $[a_0,a_m]$, implies $$\forall n\in \NN_0 \quad p=\cK^np\geq \cK^n(g\vee 0) \text{ on }[a_0,a_m],$$ yielding $$p\geq \sup_{n\in\NN_0}\cK^n(g\vee 0)=q\text{ on }[a_0,a_m].$$
\end{proof}

\subsection{Application}

\begin{ex}[Bermudan vanilla calls and puts on an asset that pays dividends at inflation rate in the Black-Scholes model] Assume $$P:=(P_t)_{t\geq 0}:=\left(\nu_{\mu t, \sigma^2t}\ast\cdot\right)_{t\geq 0},$$ where $$\sigma>0,\quad \mu:=r-\delta-\frac{\sigma^2}{2}=-\frac{\sigma^2}{2},$$ thus $P$ can be perceived as the semigroup associated to the logarithmic price process under the risk-neutral measure in the one-dimensional Black-Scholes model). We will assume that (possibly after re-scaling the time scale) $\sigma=1$. Define $$g_0:=-K+\exp, \quad g_1:K-\exp$$ (the payoffs on exercise of a one-dimensional call and a one-dimensional put option with strike price $K$, respectively). The infinitesimal generator of the Markov semigroup $P$ is $$L:f\mapsto \frac{1}{2}f'' + \mu f'=\frac{1}{2}f''+\left(r-\delta-\frac{1}{2}\right)f'.$$

Let us, for the remainder of this Example, assume that $r=\delta$. Then we have $L\exp =0$ and of course $L1=0$ ($1$ denoting the constant function $1:x\mapsto 1$). Therefore, by Lemma \ref{uniqueinterpol} the set of harmonic functions equals the space of all linear combinations of $1$ and $ \exp $, which is the set of all real multiples of the exponential function with additive constants. Thus, both $g_0$ and $g_1$ will be harmonic. 

Now suppose we are given support abscissas $a_0<\dots <a_m$ as in the previous paragraphs. 

In order to obtain the setting of Theorem \ref{fixexist} for $g=g_0$, set $h_0= a+ \exp $ for an arbitrary $a\geq -K$. Then $h_0$ will be a harmonic function dominating $g_0=\exp-K$ on $\RR$. Also, the constant function $c_0=\min_{[a_0,a_m]}\left(\exp-K\right)=e^{a_0}-K$ will be harmonic and trivially less than or equal to $g_0$ on $[a_0,a_m]$. Then the conditions of Theorem \ref{fixexist} are satisfied for $c=c_0$, $g=g_0$ and $h=h_0$, and the Theorem as well as its Corollary may be applied to establish the existence of a unique minimal fixed point for $\cK$ and its characterisation as the approximate price of the corresponding perpetual Bermudan option with exercise mesh $t$ -- the approximation being based on piecewise harmonic interpolation as outlined in the introduction to this paper. In order for $\cK$ to be financially meaningful, however, we must choose $a_0\geq \ln K$ to ensure that $c_0$, the pointwise lower bound on the elements of $Q$, is nonnegative.

Finally, in order to get the situation of Theorem \ref{fixexist} also for $g=g_1$, we now define $h_1=K-\exp$ to get a harmonic function dominating $g_1=K-\exp$ on all of $\RR$. Also, the constant function $c_1=\min_{[a_0,a_m]}\left(K-\exp \right)=K-e^{a_m}$ will be harmonic and trivially less than or equal to $g_1$ on $[a_0,a_m]$. Then the conditions of Theorem \ref{fixexist} are satisfied for the choices of $c=c_1$, $g=g_1$ and $h=h_1$, and both the Theorem and its Corollary may be applied to see that again there is a minimal fixed point for $\cK$ coinciding with the approximate price of the corresponding perpetual Bermudan option with exercise mesh $t$. In order for $\cK$ to be financially meaningful, we must this time choose $a_m\leq \ln K$ to ensure that all elements of $Q$ are pointwise nonnegative.

\end{ex}

\section{R\'eduite-based approximation of Bermudan option prices}

Suppose $P$ is a Markov semigroup on $\RR^d$ ($d\in\NN$) and $L$ is the infinitesimal generator of $P$. We will call a function $f:\RR^d\rightarrow \RR$ {\em subharmonic} if and only if $$\forall t> 0 \quad P_tf \geq f$$ holds pointwise. A function $f:\RR^d\rightarrow \RR$ will be called {\em superharmonic} if and only if $-f$ is subharmonic, and $f:\RR^d\rightarrow \RR$ will be called {\em harmomic} if it is both super- and subharmonic.

Let $\cU$ denote the operator of upper-semicontinuous regularisation, that is, for all functions $f:\RR^d\rightarrow \RR$, $$\cU f=\inf\left\{\ell \geq f \ : \ \ell:\RR^d\rightarrow\RR \text{ subharmonic}\right\}$$ (of course, this is a priori only defined as a function taking values in $\RR\cup\{--\infty\}$). Consider a harmonic function $h:\RR^d\rightarrow \RR$ and a closed (and therefore $F_{\sigma}$) set $B$ and define the {\em r\'eduite} operator $\cR=\cR_{h,B}$ on the set of all subharmonic functions $f:\RR^d\rightarrow \RR$ dominated by $h$ via $$\cR f:=\cU\left(\sup\left\{\ell \leq h \ : \ \ell:\RR^d\rightarrow\RR \text{ subharmonic}, \quad \ell \leq f \text{ on } B\right\}\right).$$ It is a well-known  result from potential theory (cf e g the work of Paul-Andr\'e Meyer \cite[Th\'eor\`eme T22]{M}) that there will be a greatest subharmonic function dominated by $f$ on $B$ and that this function will be equal to $\cR f$. Moreover, we have that $f=\cR f$ on $B$ except on a set of potential zero, in probabilistic/potential-theoretic jargon $$f=\cR f \text{ q.e. on }B,$$ where ``q.e.'' is, as usual, short-hand for ``quasi-everywhere''. Now define $$Q:=\left\{f\leq h \ : \ f:\RR^d\rightarrow\RR \text{ subharmonic}\right\}.$$ Then our definition of the r\'eduite operator $\cR$ implies $\cR f\leq h$ (as $h$ is dominating the function whose upper-semicntinuous regularisation is, according to our definition, the r\'eduite $\cR f$ of $f$) and our potential-theoretic characterisation of the r\'eduite -- as the greatest subharmonic function dominated by $f$ on $B$ -- ensures the subharmonicity of $\cR f$. Therefore, $$\cR:Q\rightarrow Q.$$ We also have that $\cU$ is monotone (in the sense that for all $f_0\leq f_1$, $\cU f_0\leq \cU f_1$) so that $\cR$ must be monotone as well (from the $\subseteq$-monotonicity of $\sup$ and the definition of $\cR$).

Hence

\begin{lem}\label{propR} Adopting the notation of the preceding paragaph, $\cR:Q\rightarrow Q$ and whenever $f_0\leq f_1$, $\cR f_0\leq \cR f_1$. 
\end{lem}

Let $g:\RR^d\rightarrow \RR$ be a subharmonic function such that $g\leq h$ and let $r>0$. The next step is going to be the consideration of the following family of operators: $$\phi_t:f\mapsto e^{-rt}P_tf\vee g$$ for $t\geq 0$. If $f\leq h$, $P_tf\leq P_t h=h$ for all $t\geq 0$, since the operators $P_t$ are positive and linear, and $h$ was assumed to be harmonic. Thus, since $g\leq h$ and $r>0$, one must have $\phi_t f\leq h$ for all $f\leq h$ and $t\geq 0$. Moreover, the operators $P_t$ preserve subharmonicity and the maximum of two subharmonic functions is subharmonic again, therefore $\phi_t f$ must be subharmonic for all subharmonic $f$. Finally, since $P_t$ is monotone, $\phi_t$ has to be monotone for all $t\geq 0$. Summarising this, we obtain

\begin{lem}\label{propphi} Using the notation introduced previously, $\phi_t:Q\rightarrow Q$ and whenever $f_0\leq f_1$, $\phi_t f_0\leq \phi_t f_1$ for all $t\geq 0$. 
\end{lem}

As a consequence, we derive from the two Lemmas \ref{propR} and \ref{propphi} the following:

\begin{cor} If we define $\cK_t:=\cR\circ\phi_t$ (adopting the notation of the previous paragraph), we have $\cK_t:Q\rightarrow Q$ and whenever $f_0\leq f_1$, $\cK_t f_0\leq \cK_t f_1$. 
\end{cor}

\begin{cor} The map $f\mapsto\cK_tf\vee 0$ is a sound iterative Bermudan option pricing algorithm for the payoff function $g\vee 0$ (in the sense of Definition \ref{algorithmclasses}).
\end{cor}

This already suffices to prove the following

\begin{Th} \label{Theorem7.2}Let $t\geq0$. Then for all $n\in\NN_0$, \begin{equation}\label{mono}{\cK_t}^{n+1}(g\vee 0)\geq {\cK_t}^{n}(g\vee 0).\end{equation} Furthermore, $$q:=\sup_{n\in \NN_0} {\cK_t}^{n}(g\vee 0)$$ (which a priori is only defined as a function with range in $\RR\cup\{+\infty\}$) is an element of $Q$ and indeed is the least nonnegative fixed point of $\cK_t$.
\end{Th}
\begin{proof}[Proof] \begin{enumerate}
\item Relation (\ref{mono}) follows from the fact that $\cK_t$ is a sound algorithm and Remark \ref{soundmonotone}.
\item Since $\cK_t$ maps $Q$ to itself, the whole sequence $\left({\cK_t}^{n}(g\vee 0)\right)_{n\in\NN_0}$ is bounded by $h$. This entails $q\leq h$ as well. Applying Beppo Levi's Theorem on swapping $\sup$ and $\int\cdot d\mu$ -- for bounded monotonely increasing sequences of measurable nonnegative functions and an arbitrary measure $\mu$ -- to the measures $P_t(\cdot, x)$, $x\in\RR^d$ and the sequence $\left({\cK_t}^{n}(g\vee 0)\right)_{n\in\NN_0}$, we can exploit the subharmonicity of the functions ${\cK_t}^{n}(g\vee 0)$, ${n\in\NN_0}$, to deduce \begin{eqnarray*}\forall x\in\RR^d\quad P_tq(x)&=&\sup_{n\in\NN_0} P_t\left({\cK_t}^{n}(g\vee 0)\right)(x)\\ &\geq& \sup_{n\in\NN_0}{\cK_t}^{n}(g\vee 0)(x) =q(x),\end{eqnarray*} which is the subharmonocity of $q$. As we have already seen, $q\leq h$, so $q\in Q$.
\item If we employ Beppo Levi's Theorem again, we can show that $\cK_t$ and $\sup_{n\in\NN_0}$ commute for bounded monotonely increasing sequences of functions. Thereby $$\cK_tq=\sup_{n\in\NN_0}\cK_t {\cK_t}^{n}(g\vee 0)=\sup_{n\in\NN}{\cK_t}^{n}(g\vee 0)=q.$$
\item That $q$ is the least nonnegative fixed point is seen as in the proof of Lemma \ref{inffixedpoint}. Any nonnegative fixed point $p$ of $\cK_t$ must be greater or equal $g\vee 0$. Therefore by the monotonicity of $\sup$ and $\cK_t$, $$\sup_{n\in\NN_0}{\cK_t}^{n}p\geq \sup_{n\in\NN_0}{\cK_t}^{n}(g\vee 0)=q.$$
\end{enumerate}
\end{proof}

\begin{ex}[Bermudan call option with equidistant exercise times in $t\cdot \NN_0$ on the weighted arithmetic average of a basket in a special Black-Scholes model] Let $\left(\beta_1,\dots,\beta_d\right)\in[0,1]$ be a convex combination and for simplicity, assume that the assets in the basket are independent and each follow the Black-Scholes model with one and the same volatility $\sigma_1=\dots=\sigma_d=:\sigma$, and let $r>0$ be the interest rate of the bond. We may assume that, possibly after a linear change of the time-scale, $\sigma=1$. Then $\left(P_t\right)_{t\geq 0}=\left(\nu_{{^t}\left(r-\frac{1}{2}\right)_{i=1}^dt,t}\ast\cdot\right)_{t\geq 0}$ is the semigroup of this Markov (even L\'evy) basket. Then one has $$L=\frac{1}{2}\Delta + \left(r-\frac{1}{2}\right)_{i=1}^d\cdot\nabla$$ (cf e g Revuz and Yor's exposition \cite{RY}), and for $$g:x\mapsto \sum_{i=1}^d \beta_i\exp\left(x_i\right)-K$$ we obtain $$Lg=\sum_{i=1}^d\left(\frac{{\beta_i}^2}{2} +\left(r-\frac{\beta_i}{2}\right)\right)\exp \left((\cdot)_i\right)$$ which is pointwise nonnegative if and only if  $$r\geq \frac{\max_{i\in\{1,\dots,d\}}{\beta_i}}{2}=\frac{1}{2d}+\frac{\sum_{i=1}^d {\beta_i}^2}{2d}. $$ Hence, if $r$ is sufficiently large, $g$ is subharmonic and we can apply the theory developed earlier in this paper, in particular Theorem \ref{Theorem7.2}.
\end{ex}

\section{Soundness and convergence rate of perpetual Bermudan option pricing via cubature}

When Nicolas Victoir studied ``asymmetric cubature formulae with few points'' \cite{V} for symmetric measures such as the Gaussian measure, the idea of (non-perpetual) Bermudan option pricing via cubature in the log-price space was born. In the following, we will discuss the soundness and convergence rate of this approach when used to price perpetual Bermudan options.

Consider a convex combination $(\alpha_1,\dots,\alpha_m)\in[0,1]^d$ (that is, $\sum_{k=1}^m \alpha_k = 1$) and $x_1,\dots,x_m\in\RR^d$. Then there is a canonical {\em weighted arithmetic average operator} $A$ associated with $\vec \alpha, \vec x$ given by $$\forall f\in\RR^\RR\quad Af=\sum_{k=1}^m\alpha_k f(\cdot-x_k).$$ Now suppose $c\in(0,1)$, $g,h:\RR^d\rightarrow \RR$, $Ag\geq g$, $Ah=h$ and $0\vee g\leq h$. Define an operator $\cD$ on the cone of nonnegative measurable functions by $$\cD:f\mapsto \left(c\cdot Af\right)\vee g.$$ 

We should stress that we explicitly allow for $g$ to take negative values.

\begin{lem} \label{discreteLemma}Adopting the previous paragraph's notation and setting $$Q:=\left\{f\leq h \ : \ Af\geq f\geq 0\right\},$$ we have the following properties of $\cD$ and $Q$:

\begin{enumerate}
\item $\cD$ is monotone (ie $\cD f\leq \cD g$ whenever $f\leq g$). 
\item $A\cD -\cD$ is nonnegative on $Q$.\label{AD-D geq 0}
\item $\cD:Q\rightarrow Q.$ \label{DmapsQtoQ} 
\item $(g\vee 0)\in Q$.

\end{enumerate}
\end{lem}
\begin{proof} \begin{enumerate}
\item Since $A$ is positive and linear, thus monotone (in the sense that for all $f_0\leq f_1$, $A f_0\leq A f_1$), it follows that $\cD$ is a composition of monotone maps, thus monotone as well. 
\item Whenever $Af\geq f$, the monotonicity of $A$ and our assumption $Ag\geq g$ imply $$A\cD f \geq A(cAf)\vee Ag\geq \left(c Af\right)\vee g =\cD f.$$ Hence $A\cD f \geq \cD f $ for all $f\in Q$. 
\item Because $h$ is nonnegative and $g\vee 0\leq h$ and due to the monotonicity of $\cD$, we have for all nonnegative $f\leq h$, $$\cD f\leq \cD h =cAh \vee g = ch\vee g= ch\vee 0\vee g\leq h,$$ thus $\cD f\leq h$ for all $f\in Q$. Also, whenever $Af\geq f\geq 0$, one gets $cAf\geq 0$ and therefore in particular $\cD f\geq 0$ for all $f\in Q$. However, we have already shown that $A\cD f \geq \cD f $ for all $f\in Q$. Summarising this, we arrive at $\cD f\in Q$ for every $f\in Q$. 
\item Because of our assumption $Ag\geq g$ and the monotonicity of $A$, we not only have $A(g\vee 0)\geq 0$, but also $A(g\vee 0)\geq g$, therefore $A(g\vee 0)\geq g\vee 0\geq 0$. However, by another assumption, $g\vee 0\leq h$. Therefore $(g\vee 0)\in Q$ as claimed.
\end{enumerate}
\end{proof}

\begin{cor} The map $f\mapsto\cD f\vee 0$ is a sound iterative Bermudan option pricing algorithm for the payoff function $g\vee 0$ (in the sense of Definition \ref{algorithmclasses}).
\end{cor}

Lemma \ref{discreteLemma} suffices to prove 

\begin{Th} \label{monoconv}For all $n\in\NN_0$, \begin{equation}\label{discretemono}\cD^{n+1}(g\vee 0)\geq \cD^n(g\vee 0)=:q_n.\end{equation} Furthermore, $$q:=\lim_{n\rightarrow \infty}\cD^n(g\vee 0)=\sup_{n\in\NN_0}\cD^n(g\vee 0)\in Q$$ and $q$ is the smallest nonnegative fixed point of $\cD$.
\end{Th}
\begin{proof}
\begin{enumerate}
\item Relation (\ref{discretemono}) follows from the soundness of $\cD$ and Remark \ref{soundmonotone}.
\item In Lemma \ref{discreteLemma}, we have not only seen that $(g\vee 0)\in Q$ and $\cD$ is closed under $Q$, but also that $A\cD f\geq \cD f$ for all $f\in Q$. Hence $$\forall n\in \NN_0 \quad A\left({\cD}^{n}(g\vee 0)\right) \geq {\cD}^{n}(g\vee 0), $$ and therefore \begin{eqnarray*}Aq &= &\sup_{n\in\NN_0} A\left({\cD}^{n}(g\vee 0)\right) \geq \sup_{n\in\NN_0}{\cD}^{n}(g\vee 0) = q,\end{eqnarray*} which means $Aq\geq q\geq 0$. 

Again because $\cD$ maps $Q$ itself and $(g\vee 0)\in Q$, we have that the whole sequence $\left({\cD}^{n}(g\vee 0)\right)_{n\in\NN_0}$ is bounded by $h$. This entails $q\leq h$ as well. 

As we have already seen, $Aq\geq q\geq 0$, so $q\in Q$.
\item Since $A$ is a weighted arithmetic average operator, $A$ and $\sup_{n\in\NN_0}$ commute in the sense that $A\left(\sup_nf_n\right)=\sup_n Af_n$ for increasing sequences of functions $(f_n)_{n\in\NN_0}$. Hence, whenever $(f_n)_{n\in\NN_0}$ is increasing, $$\cD \left(\sup_n f_n\right)= A\left(\sup_n f_n\right)\vee g= \left(\sup_n Af_n\right)\vee g =  \sup_n \left(Af_n\vee g\right)=\sup_n\cD f_n,$$ i e $\cD$ and $\sup$ commute for bounded monotonely increasing sequences of functions. Thereby $$\cD q=\sup_{n\in\NN_0}\cD{\cD}^{n}(g\vee 0)=\sup_{n\in\NN} {\cD}^{n}(g\vee 0)=q.$$
\item Just as in the proof of Lemma \ref{inffixedpoint}, we see that $q$ is the minimal nonnegative fixed point. For, any nonnegative fixed point $p$ of $\cD$ must be greater or equal $g\vee 0$. Thus, by the monotonicity of $\sup$ and $\cD$, $$\sup_{n\in\NN_0}{\cD}^{n}p\geq \sup_{n\in\NN_0}{\cD}^{n}(g\vee 0)=q.$$
\end{enumerate}
\end{proof}

\begin{lem} \label{gLemma} Using the previous Theorem's notation, we have for all $x\in\RR^d$ and $n\in\NN_0$, if $q_{n+1}(x)=g(x)$, then $q_n(x)=g(x)$.
\end{lem}
\begin{proof} By the monotonicity of the sequence $(q_n)_{n\in\NN_0}$ (Theorem \ref{monoconv}), we have $$g(x)\leq q_0(x)\leq q_n(x)\leq q_{n+1}(x).$$
\end{proof}

\begin{Th} For all $n\in\NN$, $$\left\|q_{n+1}-q_n\right\|_{C^0\left(\RR^d,\RR\right)}\leq c\cdot \left\|q_{n}-q_{n-1}\right\|_{C^0\left(\RR^d,\RR\right)}.$$
\end{Th}

\begin{proof} The preceding Lemma \ref{gLemma} yields \begin{eqnarray*} \left\|q_{n+1}-q_n\right\|_{C^0\left(\RR^d,\RR\right)}&=& \left\|q_{n+1}-q_n\right\|_{C^0\left(\left\{c\cdot Aq_{n}>g\right\},\RR\right)}\\ &=&\left\|c\cdot Aq_{n}-\left(\left(c\cdot A q_{n-1}\right)\vee g\right)\right\|_{C^0\left(\left\{c\cdot Aq_{n}>g\right\},\RR\right)}\end{eqnarray*} via the definition of $q_{i+1}$ as $\left(cAq_i\right)\vee g$ for $i=n$ and $i=n+1$. But the last equality implies \begin{eqnarray*}\left\|q_{n+1}-q_n\right\|_{C^0\left(\RR^d,\RR\right)}&\leq&\left\|c\cdot Aq_{n}-c\cdot A q_{n-1}\right\|_{C^0\left(\left\{c\cdot Aq_{n}>g\right\},\RR\right)}\\ &\leq&\left\|c\cdot Aq_{n}-c\cdot A q_{n-1}\right\|_{C^0\left(\RR^d,\RR\right)}.\end{eqnarray*} Since $A$ is linear as well as an $L^\infty$-contraction (and therefore a $C^0$-contraction, too), we finally obtain $$ \left\|q_{n+1}-q_n\right\|_{C^0\left(\RR^d,\RR\right)}\leq c\left\|A\left(q_{n}-q_{n-1}\right)\right\|_{C^0\left(\RR^d,\RR\right)}\leq c\left\|q_{n}-q_{n-1}\right\|_{C^0\left(\RR^d,\RR\right)}.$$
\end{proof}

\begin{ex}[Bermudan put option with equidistant exercise times in $t\cdot \NN_0$ on the weighted arithmetic average of a basket in a discrete Markov model with a discount factor $c=e^{-rt}$ for $r>0$] Let $\beta_1,\dots,\beta_d\in[0,1]$ be a convex combination and assume that $A$ is such that \begin{equation}\label{discretecondition}\forall i\in \{1,\dots,d\}\quad \sum_{k=1}^m\alpha_ke^{-(x_k)_i}=1,\end{equation} then the functions $$g:x\mapsto K-\sum_{i=1}^d\beta_i \exp\left(x_i\right)$$ and $h:=K$ (where $K\geq 0$) satisfy the equations $Ah=h$ and $Ag=g$, respectively. Moreover, by definition $g\leq h$, thus $0\vee g\leq h$. Then we know that the (perpetual) Bermudan option pricing algorithm that iteratively applies $\cD$ to the payoff function $g\vee 0$ on the $\log$-price space, will increase monotonely and will have a limit which is the smallest nonnegative fixed point of $\cD$. Moreover, the convergence is linear and the contraction rate can be bounded by $c$.

The condition (\ref{discretecondition}) can be achieved by a change of the time scale (which ultimately leads to different cubature points for the distribution of the asset price)
\end{ex}

One might also be interested in determining the convergence rate for the approximation of non-perpetual American option prices based on non-perpetual Bermudan option pricing via cubature. After proving a series of Lemmas we will end up with a Theorem that asserts linear convergence and also provides bounds for the convergence factor. 

From now on, $c$ and $A$ will no longer be fixed but their r\^ole will be played by $e^{-rt}$ and $P_t$ (for $t\in s\NN_0$ where $s>0$ shall be fixed) respectively, where $r>0$ and $\left(P_{s\cdot m}\right)_{m\in\NN_0}$ describes a Markov chain on $\RR^d$ (By the Chapman-Komogorov equation this is tantamount to $\forall s,t\geq 0\quad P_sP_t=P_{s+t}$).

{\bf Acknowledgements.} The author would like to thank the German Academic Exchange Service for the pre-doctoral research grant he received ({\em Doktorandenstipendium des Deutschen Akademischen Austauschdienstes}) and the German National Academic Foundation ({\em Studienstiftung des deutschen Volkes}) for their generous support in both financial and non-material terms. 

Moreover, he owes a huge debt of gratitude to his supervisor, Professor Terry J Lyons, for numerous extremely helpful discussions, as well as to Dr Ben Hambly and Professor Alexander Schied for their constructive comments on a previous version of this paper.

\end{document}